# Isogeometric finite element analysis of functionally graded plates using a refined plate theory

**H**. **Nguyen-Xuan**[1,2*], **Loc V**. **Tran**[3], **Chien H**. **Thai**[3], **S**. **Kulasegaram**[1], **S.P.A. Bordas**[1]

[1]Institute of Mechanics and Advanced Materials, School of Engineering, Cardiff University, Queen's Buildings, The Parade, Cardiff CF24 3AA, UK

[2]Department of Mechanics, Faculty of Mathematics & Computer Science, University of Science Ho Chi Minh City, Vietnam

[3]Division of Computational Mechanics, Ton Duc Thang University Ho Chi Minh City, Vietnam

**Abstract**

We propose in this paper a novel inverse tangent transverse shear deformation formulation for functionally graded material (FGM) plates. The isogeometric finite element analysis (IGA) of static, free vibration and buckling problems of FGM plates is then addressed using a refined plate theory (RPT). The RPT enables us to describe the non-linear distribution of shear stresses through the plate thickness without any requirement of shear correction factors (SCF). IGA utilizes basis functions, namely B-splines or non-uniform rational B-splines (NURBS), which achieve easily the smoothness of any arbitrary order. It hence satisfies the $C^1$ requirement of the RPT model. The present method approximates the displacement field of four degrees of freedom per each control point and retains the computational efficiency while ensuring the high accuracy in solution.

## 1. Introduction

FGM-a mixture of two distinct material phases: ceramic and metal, e.g. Figure 1, are very well capable of reducing thermal stresses, resisting high temperature environment and corrosion coatings. FG plates are therefore suitable for applications in aerospace structures, nuclear plants and semiconductor technologies. They have hence been received great attention by many researchers.

---

[*]Corresponding author. *Email address*: nxhung@hcmus.edu.vn (H. Nguyen-Xuan)

Besides the extensive application of FG plates in engineering structures, a lot of plate theories have been developed to analyze the thermo-mechanical behavior of such structures. The Classical Plate Theory (CPT) relied on the Kirchoff-Love assumptions merely to provide acceptable results for the thin plate. The First Order Shear Deformation Theory (FSDT) based on Reissner [1] and Mindlin [2], which takes into account the shear effect, was therefore developed. In FSDT-based finite element models, it is necessary to use some improved techniques such as reduced integration (RI) [3], mixed interpolation of tensorial components (MITC) [4], Mindlin-type plate element with improved transverse shear(MIN) [5], discrete shear gap (DSG) [6-8] elements, etc, to overcome the shear locking phenomenon. In addition, with the linear in-plane displacement assumption through plate thickness, shear strain/stress obtained from FSDT distributes inaccurately and does not satisfy the traction free boundary conditions at the plate surfaces. The shear correction factors (SCF) are hence required to rectify the unrealistic shear strain energy part. The values of SCF are quite dispersed through each problem and may be difficult to determine [9]. To ensure the curved distribution of shear stress, various kinds of Higher-Order Shear Deformable Theory (HSDT), which includes higher-order terms in the approximation of the displacement field, have then been devised such as Third-Order Shear Deformation Theory (TSDT) [10-14], trigonometric shear deformation theory [15-19], exponential shear deformation theory [20-22], refined plate theory (RPT) based on two unknown functions of transverse deflection [23,24] and so on. The two variable RPT model was first proposed by Shimpi [25] using only two unknown variables for the isotropic plate and then extended for the orthotropic plate [26,27]. It is worth mentioning that the HSDT models provide better results and yield more accurate and stable solutions (e.g. inter-laminar stresses and displacements) [28,29] than the FSDT ones without requirement the SCF. However, the HSDT requires the $C^1$-continuity of generalized displacement field leading to the second-order derivative of the stiffness formulation and it causes the obstacles in the standard finite element formulations. Several $C^0$ continuous elements [30-34] were proposed or Hermite interpolation function with the $C^1$-continuity was added for just specific approximation of transverse displacement [10]. It may produce extra unknown variables including derivative of deflection $w_{,x}$, $w_{,y}$, $w_{,xy}$ [35] leading to increase in the computational cost. In this paper, we show that $C^1$-continuous elements will be easily achieved by IGA without any additional variables.

Isogeometric approach (IGA) has been recently proposed by Hughes et al. [36] to closely link the gap between Computer Aided Design (CAD) and Finite Element Analysis (FEA). The basic idea is that the IGA uses the same non-uniform rational B-Spline (NURBS) functions in describing the exact geometry of problem and constructing finite approximation for analysis. It is well known that NURBS functions provide a flexible way to make refinement, de-refinement, and degree elevation [37]. They enable us to easily achieve the smoothness of arbitrary continuity order in comparison with the traditional FEM. Hence, IGA naturally verifies the $C^1$-continuity of plates based on the HSDT assumption, which is interested in this study. The IGA has been well known and widely applied to various practical problems [38-44] and so on.

In this paper, a novel inverse tangent transverse shear deformation plate theory with four parameters of displacement field are proposed to study the behavior of the FGM plates based on NURBS-based IGA approximation. The RPTs are independent on SCF and free from shear locking. The material property changing continuously through plate thickness is homogenized by the rule of mixture and the Mori-Tanaka homogenization technique. IGA utilizes the NURBS functions which achieve easily the smoothness with arbitrary continuous order. It helps the present method to naturally satisfy the $C^1$ continuous requirement of this plate theory. Numerous numerical examples are provided to illustrate the effectiveness of the present formulation in comparison with other published models.

The paper is outlined as follows. Next section introduces the novel RPT for FGM plates. In section 3, the formulation of plate theory based on IGA is described. The numerical results and discussions are provided in section 4. Finally, this article is closed with some concluding remarks.

## 2. The novel refined plate theory for FGM plates

### 2.1. Functionally graded material

Functionally graded material is a composite material which is created by mixing two distinct material phases which are often ceramic at the top and metal at the bottom. In this paper, two homogenous models have been used to estimate the effective properties of the FGM including the rule of mixture and the Mori-Tanaka technique. The volume fraction of the ceramic and metal phase is assumed continuous variation through thickness as the following exponent function [10]

$$V_c(z) = \left(\frac{1}{2} + \frac{z}{h}\right)^n, \quad V_m = 1 - V_c \qquad (1)$$

where subscripts $m$ and $c$ refer to the metal and ceramic constituents, respectively. Eq. (1) denotes that the volume fraction varies through the thickness based on the power index $n$.

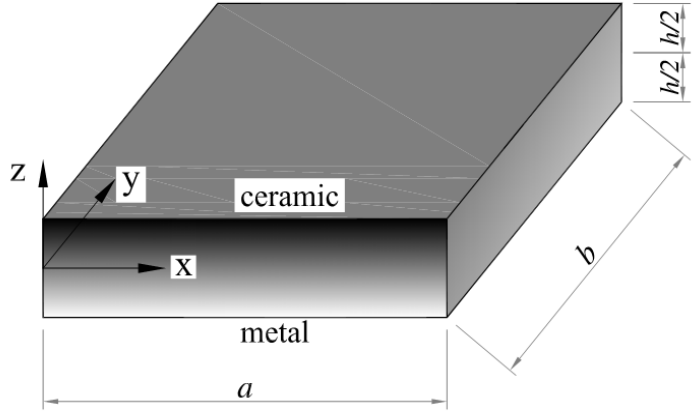

Figure 1: The functionally graded plate model.

Then, the effective properties of material according to the rule of mixture [10] are given by

$$P_e = P_c V_c(z) + P_m V_m(z) \tag{2}$$

where $P_c, P_m$ denote the individual material's properties of the ceramic and the metal, respectively including the Young's modulus $E$, Poisson's ratio $\nu$, density $\rho$.

However, the rule of mixture does not consider the interactions among the constituents [46,47]. So, the Mori-Tanaka technique [45] is then developed to take into account these interactions with the effective bulk and shear modulus defined from following relations:

$$\begin{aligned} \frac{K_e - K_m}{K_c - K_m} &= \frac{V_c}{1 + V_m \frac{K_c - K_m}{K_m + 4/3\mu_m}} \\ \frac{\mu_e - \mu_m}{\mu_c - \mu_m} &= \frac{V_c}{1 + V_m \frac{\mu_c - \mu_m}{\mu_m + f_1}} \end{aligned} \tag{3}$$

where $f_1 = \dfrac{\mu_m(9K_m + 8\mu_m)}{6(K_m + 2\mu_m)}$. Then, the effective values of Young's modulus $E_e$ and Poisson's ratio $\nu_e$ are given by:

$$E_e = \frac{9K_e \mu_e}{3K_e + \mu_e}, \qquad \nu_e = \frac{3K_e - 2\mu_e}{2(3K_e + \mu_e)} \tag{4}$$

Figure 2 illustrates the comparison of the effective Young's modulus of Al/ZrO$_2$ FGM plate calculated by the rule of mixture and the Mori-Tanaka scheme via the power index $n$. It is noted that with homogeneous materials, two models give the same values, however as material becomes inhomogeneous, the effective property through the thickness of the former is higher than the latter one.

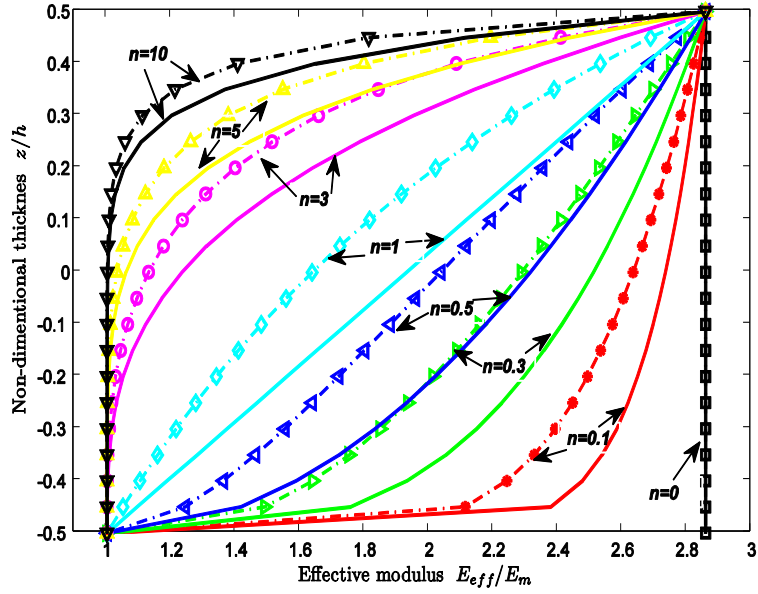

Figure 2. The effective modulus of Al/ZrO$_2$ FGM plate computed by the rule of mixture (in solid line) and the Mori-Tanaka (in dash dot line).

*2.2. The novel inverse tangent transverse shear deformation plate theories* 

A general displacement field based on refined plate theory according to Shimpi [25] is defined as follows

$$\begin{aligned}
u(x,y,z) &= u_0 - zw_{b,x} + g(z)w_{s,x} \\
v(x,y,z) &= v_0 - zw_{b,y} + g(z)w_{s,y} \quad , \quad \left(\frac{-h}{2} \le z \le \frac{h}{2}\right) \\
w(x,y) &= w_b + w_s
\end{aligned} \quad (5)$$

where $u_0, v_0$ are the membrane displacements along the *x, y* axis respectively and $w_b, w_s$ are the bending and shear components of deflection and the function $g(z) = f(z) - z$ is used to describe the distribution of the transverse shear strains and stresses through the plate thickness. The function $f(z)$ is chosen to satisfy the tangential zero value at $z = \pm h/2$. It means that traction-free boundary condition at the top and bottom plate surfaces is automatically satisfied. In this study, beyond some existing shape functions [10,20,21,48], we propose two inverse tangent formulations for $f(z)$ as listed in Table 1.

Table 1: The various forms of shape function.

| Model | $f(z)$ | $f'(z)$ |
|---|---|---|
| Reddy [10] | $z - \frac{4}{3}z^3/h^2$ | $1 - 4z^2/h^2$ |
| Karama [20] | $ze^{-2(z/h)^2}$ | $(1 - \frac{4}{h^2}z^2)e^{-2(z/h)^2}$ |
| Arya [21] | $\sin(\frac{\pi}{h}z)$ | $\frac{\pi}{h}\cos(\frac{\pi}{h}z)$ |
| Nguyen-Xuan [48] | $\frac{7}{8}z - \frac{2}{h^2}z^3 + \frac{2}{h^4}z^5$ | $\frac{7}{8} - \frac{6}{h^2}z^2 + \frac{10}{h^4}z^4$ |
| Proposed model 1 | $h\arctan(\frac{2}{h}z) - z$ | $(1 - (\frac{2}{h}z)^2)/(1 + (\frac{2}{h}z)^2)$ |
| Proposed model 2 | $\arctan(\sin(\frac{\pi}{h}z))$ | $\frac{\pi}{h}\cos(\frac{\pi}{h}z)/(1 + \sin^2(\frac{\pi}{h}z))$ |

The relationship between strains and displacements is described by

$$\boldsymbol{\varepsilon} = [\varepsilon_{xx}\ \varepsilon_{yy}\ \gamma_{xy}]^T = \boldsymbol{\varepsilon}_0 + z\boldsymbol{\kappa}_b + g(z)\boldsymbol{\kappa}_s$$

$$\boldsymbol{\gamma} = [\gamma_{xz}\ \gamma_{yz}]^T = f'(z)\boldsymbol{\varepsilon}_s \tag{6}$$

where

$$\boldsymbol{\varepsilon}_0 = \begin{bmatrix} u_{0,x} \\ v_{0,y} \\ u_{0,y} + v_{0,x} \end{bmatrix}, \boldsymbol{\kappa}_b = -\begin{bmatrix} w_{b,xx} \\ w_{b,yy} \\ 2w_{b,xy} \end{bmatrix}, \boldsymbol{\kappa}_s = \begin{bmatrix} w_{s,xx} \\ w_{s,yy} \\ 2w_{s,xy} \end{bmatrix}, \boldsymbol{\varepsilon}_s = \begin{bmatrix} w_{s,x} \\ w_{s,y} \end{bmatrix} \tag{7}$$

A weak form of the static model for the plates under transverse loading $q_0$ can be briefly expressed as:

$$\int_\Omega \delta\boldsymbol{\varepsilon}^T \mathbf{D}^b \boldsymbol{\varepsilon}\,d\Omega + \int_\Omega \delta\boldsymbol{\gamma}^T \mathbf{D}^s \boldsymbol{\gamma}\,d\Omega = \int_\Omega \delta w q_0\,d\Omega \tag{8}$$

where

$$\mathbf{D}^b = \begin{bmatrix} \mathbf{A} & \mathbf{B} & \mathbf{E} \\ \mathbf{B} & \mathbf{D} & \mathbf{F} \\ \mathbf{E} & \mathbf{F} & \mathbf{H} \end{bmatrix}$$

$$A_{ij}, B_{ij}, D_{ij}, E_{ij}, F_{ij}, H_{ij} = \int_{-h/2}^{h/2} (1, z, z^2, g(z), zg(z), g^2(z))Q_{ij}\,dz \tag{9}$$

$$D_{ij}^s = \int_{-h/2}^{h/2} [f'(z)]^2 G_{ij}\,dz$$

and the material matrices are given as

$$\mathbf{Q} = \frac{E_e}{1-v_e^2}\begin{bmatrix} 1 & v_e & 0 \\ v_e & 1 & 0 \\ 0 & 0 & (1-v_e)/2 \end{bmatrix} \qquad (10)$$

$$\mathbf{G} = \frac{E_e}{2(1+v_e)}\begin{bmatrix} 1 & 0 \\ 0 & 1 \end{bmatrix}$$

For the free vibration analysis of the plates, weak form can be derived from the following dynamic equation

$$\int_\Omega \delta\boldsymbol{\varepsilon}^T \mathbf{D}^b \boldsymbol{\varepsilon} d\Omega + \int_\Omega \delta\boldsymbol{\gamma}^T \mathbf{D}^s \boldsymbol{\gamma} d\Omega = \int_\Omega \delta\tilde{\mathbf{u}}^T \mathbf{m} \ddot{\tilde{\mathbf{u}}} d\Omega \qquad (11)$$

where **m** - the mass matrix is calculated according to the consistent form

$$\mathbf{m} = \begin{bmatrix} \mathbf{I_0} & 0 & 0 \\ 0 & \mathbf{I_0} & 0 \\ 0 & 0 & \mathbf{I_0} \end{bmatrix} \text{ where } \mathbf{I_0} = \begin{bmatrix} I_1 & I_2 & I_4 \\ I_2 & I_3 & I_5 \\ I_4 & I_5 & I_6 \end{bmatrix} \qquad (12)$$

$$(I_1, I_2, I_3, I_4, I_5, I_6) = \int_{-h/2}^{h/2} \rho(z)\left(1, z, z^2, g(z), zg(z), g^2(z)\right) dz. \qquad (13)$$

and

$$\tilde{\mathbf{u}} = \begin{Bmatrix} \mathbf{u}_1 \\ \mathbf{u}_2 \\ \mathbf{u}_3 \end{Bmatrix}, \quad \mathbf{u}_1 = \begin{Bmatrix} u_0 \\ -w_{b,x} \\ w_{s,x} \end{Bmatrix}; \mathbf{u}_2 = \begin{Bmatrix} v_0 \\ -w_{b,y} \\ w_{s,y} \end{Bmatrix}; \mathbf{u}_3 = \begin{Bmatrix} w \\ 0 \\ 0 \end{Bmatrix} \qquad (14)$$

For the buckling analysis, a weak form of the plate under the in-plane forces can be expressed as:

$$\int_\Omega \delta\boldsymbol{\varepsilon}^T \mathbf{D}^b \boldsymbol{\varepsilon} d\Omega + \int_\Omega \delta\boldsymbol{\gamma}^T \mathbf{D}^s \boldsymbol{\gamma} d\Omega + \int_\Omega \nabla^T \delta w \mathbf{N}_0 \nabla w d\Omega = 0 \qquad (15)$$

where $\nabla^T = [\partial/\partial x \ \partial/\partial y]^T$ is the gradient operator and $\mathbf{N}_0 = \begin{bmatrix} N_x^0 & N_{xy}^0 \\ N_{xy}^0 & N_y^0 \end{bmatrix}$ is a matrix related to the pre-buckling loads.

## 3. The FGM plate formulation based on NURBS basis functions

*3.1. A brief of NURBS functions*

A knot vector $\Xi = \{\xi_1, \xi_2, ..., \xi_{n+p+1}\}$ is defined as a sequence of knot value $\xi_i \in R$, $i = 1,...n+p$. If the first and the last knots are repeated $p+1$ times, the knot vector is called open knot. A B-spline basis function is $C^\infty$ continuous inside a knot span and $C^{p-1}$ continuous at a single knot. Thus, as $p \geq 2$ the present approach always satisfies $C^1$-requirement in approximate formulations of RPT.

The B-spline basis functions $N_{i,p}(\xi)$ are defined by the following recursion formula

$$N_{i,p}(\xi) = \frac{\xi - \xi_i}{\xi_{i+p} - \xi_i} N_{i,p-1}(\xi) + \frac{\xi_{i+p+1} - \xi}{\xi_{i+p+1} - \xi_{i+1}} N_{i+1,p-1}(\xi)$$

$$\text{as } p = 0, \; N_{i,0}(\xi) = \begin{cases} 1 & if \; \xi_i < \xi < \xi_{i+1} \\ 0 & \text{otherwise} \end{cases}$$

(16)

By the tensor product of basis functions in two parametric dimensions $\xi$ and $\eta$ with two knot vectors $\Xi = \{\xi_1, \xi_2, ..., \xi_{n+p+1}\}$ and $\mathbf{H} = \{\eta_1, \eta_2, ..., \eta_{m+q+1}\}$, the two-dimensional B-spline basis functions are obtained

$$N_A(\xi, \eta) = N_{i,p}(\xi) M_{j,q}(\eta) \tag{17}$$

Figure 3 illustrates the set of one-dimensional and two-dimensional B-spline basis functions according to open uniform knot vector $\Xi = \{0,0,0,0,\frac{1}{2},1,1,1,1\}$.

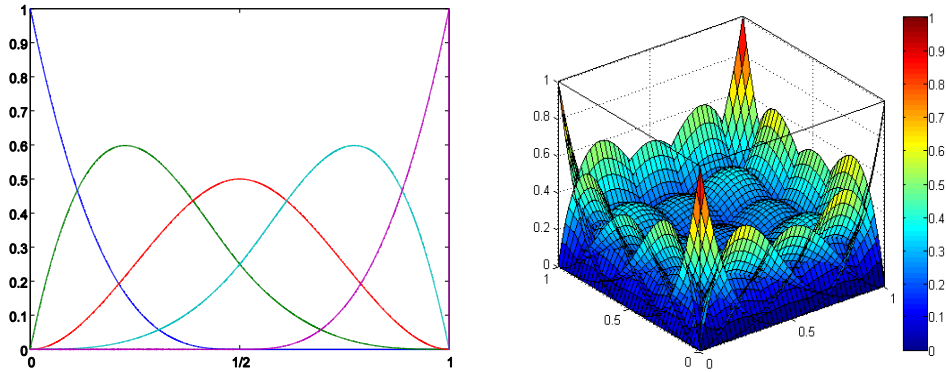

Figure 3. 1D and 2D B-spline basis functions.

To present exactly some curved geometry (e.g. circles, cylinders, spheres, etc.) the non-uniform rational B-splines (NURBS) functions are used. Be different from B-spline, each control point of NURBS has additional value called an individual weight $w_A$ [36]. Then the NURBS functions can be expressed as

$$R_A(\xi,\eta) = \frac{N_A w_A}{\sum_A^{m \times n} N_A(\xi,\eta) w_A} \tag{18}$$

It can be noted that B-spline function is only the special case of the NURBS function when the individual weight of control point is constant.

*3.2. A novel RPT formulation based on NURBS approximation* 

Using the NURBS basis functions above, the displacement field **u** of the plate is approximated as

$$\mathbf{u}^h(\xi,\eta) = \sum_A^{m \times n} R_A(\xi,\eta) \mathbf{q}_A \tag{19}$$

where $\mathbf{q}_A = [u_{0A} \ v_{0A} \ w_{bA} \ w_{sA}]^T$ is the vector of nodal degrees of freedom associated with the control point A.

Substituting Eq. (19) into Eq. (7), the in-plane and shear strains can be rewritten as:

$$\left[\boldsymbol{\varepsilon}_0^T \ \boldsymbol{\kappa}_b^T \ \boldsymbol{\kappa}_s^T \ \boldsymbol{\varepsilon}_s^T\right]^T = \sum_{A=1}^{m \times n} \left[\left(\mathbf{B}_A^m\right)^T \ \left(\mathbf{B}_A^{b1}\right)^T \ \left(\mathbf{B}_A^{b2}\right)^T \ \left(\mathbf{B}_A^s\right)^T\right]^T \mathbf{q}_A \tag{20}$$

in which

$$\mathbf{B}_A^m = \begin{bmatrix} R_{A,x} & 0 & 0 & 0 \\ 0 & R_{A,y} & 0 & 0 \\ R_{A,y} & R_{A,x} & 0 & 0 \end{bmatrix}, \quad \mathbf{B}_A^{b1} = -\begin{bmatrix} 0 & 0 & R_{A,xx} & 0 \\ 0 & 0 & R_{A,yy} & 0 \\ 0 & 0 & 2R_{A,xy} & 0 \end{bmatrix},$$

$$\mathbf{B}_A^{b2} = \begin{bmatrix} 0 & 0 & 0 & R_{A,xx} \\ 0 & 0 & 0 & R_{A,yy} \\ 0 & 0 & 0 & 2R_{A,xy} \end{bmatrix}, \quad \mathbf{B}_A^s = \begin{bmatrix} 0 & 0 & 0 & R_{A,x} \\ 0 & 0 & 0 & R_{A,y} \end{bmatrix} \tag{21}$$

Substituting Eq. (20) into Eqs.(8), (11) and (15), the formulations of static, free vibration and buckling problem are rewritten in the following form

$$\mathbf{Kq} = \mathbf{F} \tag{22}$$

$$\left(\mathbf{K} - \omega^2 \mathbf{M}\right)\mathbf{q} = \mathbf{0} \tag{23}$$

$$\left(\mathbf{K} - \lambda_{cr} \mathbf{K}_g\right)\mathbf{q} = \mathbf{0} \tag{24}$$

where the global stiffness matrix **K** is given by

$$\mathbf{K} = \int_\Omega \begin{Bmatrix} \mathbf{B}^m \\ \mathbf{B}^{b1} \\ \mathbf{B}^{b2} \end{Bmatrix}^T \begin{bmatrix} \mathbf{A} & \mathbf{B} & \mathbf{E} \\ \mathbf{B} & \mathbf{D} & \mathbf{F} \\ \mathbf{E} & \mathbf{F} & \mathbf{H} \end{bmatrix} \begin{Bmatrix} \mathbf{B}^m \\ \mathbf{B}^{b1} \\ \mathbf{B}^{b2} \end{Bmatrix} + \mathbf{B}^{sT} \mathbf{D}^s \mathbf{B}^s \mathrm{d}\Omega \qquad (25)$$

and the load vector is computed by

$$\mathbf{F} = \int_\Omega q_0 \mathbf{R} \mathrm{d}\Omega \qquad (26)$$

where

$$\mathbf{R} = \begin{bmatrix} 0 & 0 & R_A & R_A \end{bmatrix} \qquad (27)$$

the global mass matrix **M** is expressed as

$$\mathbf{M} = \int_\Omega \tilde{\mathbf{R}}^T \mathbf{m} \tilde{\mathbf{R}} \mathrm{d}\Omega \qquad (28)$$

where

$$\tilde{\mathbf{R}} = \begin{Bmatrix} \mathbf{R}_1 \\ \mathbf{R}_2 \\ \mathbf{R}_3 \end{Bmatrix}, \quad \mathbf{R}_1 = \begin{bmatrix} R_A & 0 & 0 & 0 \\ 0 & 0 & -R_{A,x} & 0 \\ 0 & 0 & 0 & R_{A,x} \end{bmatrix};$$

$$\mathbf{R}_2 = \begin{bmatrix} 0 & R_A & 0 & 0 \\ 0 & 0 & -R_{A,y} & 0 \\ 0 & 0 & 0 & R_{A,y} \end{bmatrix}; \mathbf{R}_3 = \begin{bmatrix} 0 & 0 & R_A & R_A \\ 0 & 0 & 0 & 0 \\ 0 & 0 & 0 & 0 \end{bmatrix} \qquad (29)$$

the geometric stiffness matrix is

$$\mathbf{K}_g = \int_\Omega (\mathbf{B}^g)^T \mathbf{N}_0 \mathbf{B}^g \mathrm{d}\Omega \qquad (30)$$

where

$$\mathbf{B}_A^g = \begin{bmatrix} 0 & 0 & R_{A,x} & R_{A,x} \\ 0 & 0 & R_{A,y} & R_{A,y} \end{bmatrix} \qquad (31)$$

in which $\omega$, $\lambda_{cr} \in R^+$ are the natural frequency and the critical buckling value, respectively.

It is observed from Eq. (25) that the SCF is no longer required in the stiffness formulation. Herein, $\mathbf{B}_A^{b1}$ and $\mathbf{B}_A^{b2}$ contain the second-order derivative of the shape function. Hence, it requires $C^1$- continuous element in approximate formulations. It is now interesting to note that our present formulation based on IGA naturally satisfies $C^1$-continuity from the theoretical/mechanical viewpoint of FGM plates [43, 29]. In our work, the basis functions are $C^{p-1}$ continuous. Therefore, as $p \geq 2$, the present approach always satisfies $C^1$-requirement in approximate formulations based on the proposed RPT.

## 4. Results and discussions

In this section, the plates with two kinds of shape such as square and circle are modeled. The FGM plates are made from Aluminum/Alumina (Al/Al$_2$O$_3$) or Aluminum/Zirconia (Al/ZrO$_2$) and the properties of which are listed in Table 2. Numerical results are obtained by IGA with full $(p+1) \times (q+1)$ Gauss points. In addition, two types of boundary condition are applied including:

Simply supported (S):

$$v_0 = w_b = w_s = 0 \quad \text{at} \quad x = 0, a$$
$$u_0 = w_b = w_s = 0 \quad \text{at} \quad y = 0, b \tag{32}$$

Clamped (C):

$$u_0 = v_0 = w_b = w_s = w_{b,n} = w_{s,n} = 0 \tag{33}$$

The Dirichlet BCs on $u_0, v_0, w_b$ and $w_s$ is easily treated as in the standard FEM. However, for the derivatives $w_{b,n}, w_{s,n}$ the enforcement of Dirichlet BCs can be solved in a simple and effective way [49]. The idea is as follows. The derivatives can be included in a compact form of the normal slope at the boundary:

$$\frac{\partial w}{\partial n} = \lim_{\Delta n \to 0} \frac{w(n(\mathbf{x}_C) + \Delta n) - w(n(\mathbf{x}_C))}{\Delta n} = 0 \tag{34}$$

As $w(n(\mathbf{x}_C)) = 0$ according to Eq. (33), Eq. (34) leads to impose the same boundary values, i.e, zero values, on the deflection variable at control points $\mathbf{x}_A$ which is adjacent to the boundary control points $\mathbf{x}_C$. It can be observed that, implementing the essential boundary condition using this method is very simple in IGA compare to other numerical methods.

Table 2: Material properties.

|  | Al | SiC | ZrO$_2$-1 | ZrO$_2$-2 | Al$_2$O$_3$ |
|---|---|---|---|---|---|
| $E$ (GPa) | 70 | 427 | 200 | 151 | 380 |
| $\nu$ | 0.3 | 0.17 | 0.3 | 0.3 | 0.3 |
| $\rho$ (kg/m$^3$) | 2707 | - | 5700 | 3000 | 3800 |

*4.1 Convergence study*

Let us consider the simply supported Al/SiC square FGM plate shown in Figure 4a, for which properties are given in

Table 2. The plate is subjected to a sinusoidal pressure defined as $q_0 \sin(\frac{\pi x}{a})\sin(\frac{\pi y}{a})$ at the top surface.

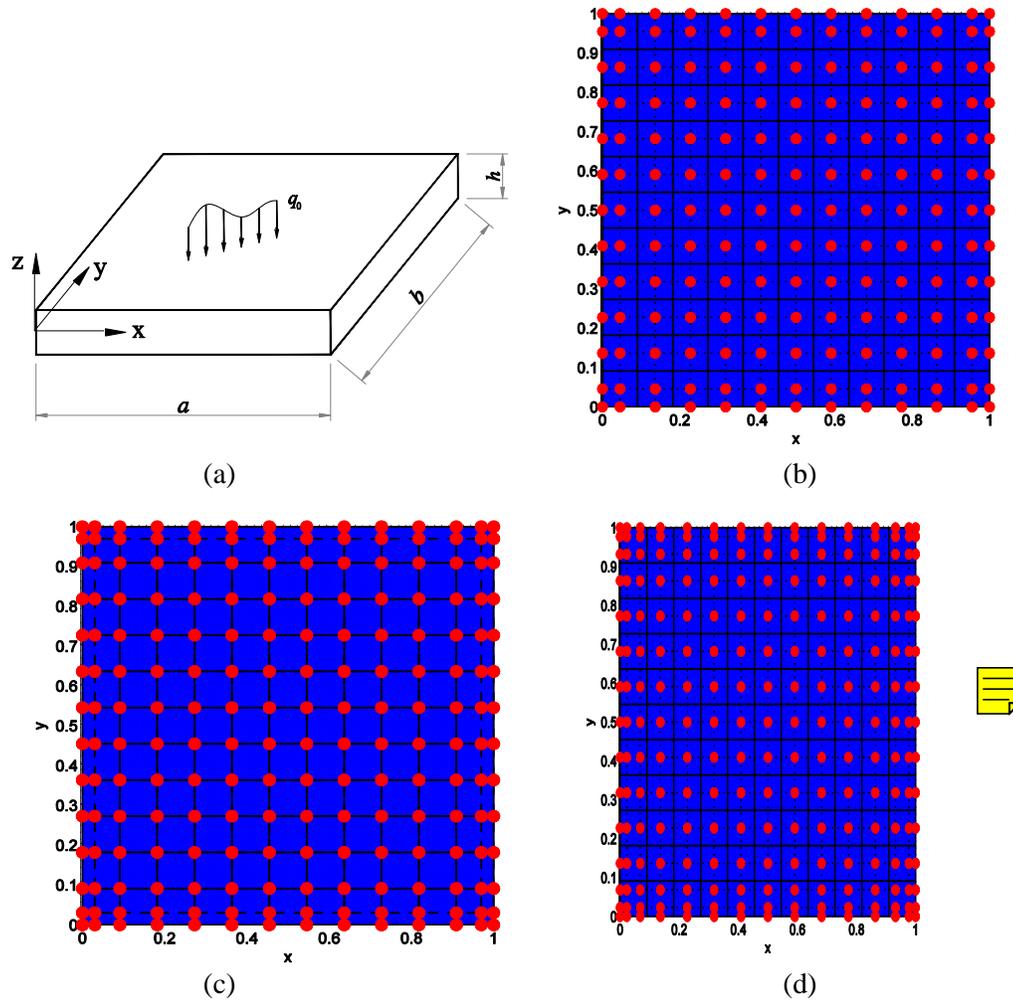

Figure 4. Square plate: (a) The plate geometry; (b), (c), (d): meshing of 11x11 quadratic, cubic and quartic elements, respectively.

A convergent study of transverse displacement by quadratic ($p = 2$) cubic ($p = 3$) and quartic ($p = 4$) elements is depicted in Figure 6 according to $n = 1$ and 6, respectively. It is observed that, as number of element increases the obtained results converge to exact solutions from 3D deformation model by Vel and Batra [46]. IGA, moreover, gains the super-convergence with the discrepancy between meshing of 5x5 and 25x25 around 0.05% as $p \geq 3$. Here, using the RPT with function $f(z)$ proposed by Reddy [10], the present method, which just uses 11x11 cubic NURBS elements shown in Figure 4c, produces an ultra-accurate solution that is very close to the exact solution with very small error around 0.02%. Therefore, in the next problems, the meshing of 11x11 cubic NURBS elements is used.

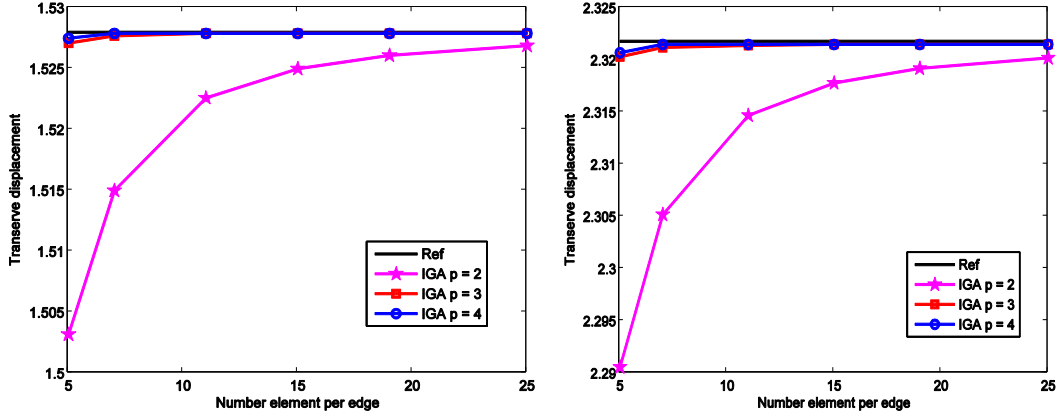

Figure 5. Comparison of present result with analytical solution of Vel and Batra according to power index $n = 1$ and 6, respectively.

Next, the shear locking phenomenon is investigated for an isotropic plate which is subjected to a uniform transverse load $q_0$. Figure 6 reveals the relation between the central deflection and the length to thickness ratio under full simply supported and full clamped conditions. By using various $f(z)$ functions in [21, 48] and the two proposed models 1 & 2, it is seen that all obtained results are the same for this problem. Those solutions are in good agreement with that using TSDT based on the Mindlin plate model with 5DOFs/node provided in [29] for moderate thick and thick plates. However, as plate becomes very thin ($a/h>1000$) the obtained results from TSDT with 5DOFs are not asymptotic to that of CPT [51]. It is called shear-locking phenomenon. Herein, with two approximated variables for transverse displacement according to Eq. (5), the shear strains/stresses are obtained independently on the bending component. As a result, the RPT has strong similarity with the CPT and will be free of shear-locking. The present results therefore match well with the thin plate result [51], even at $a/h = 10^6$.

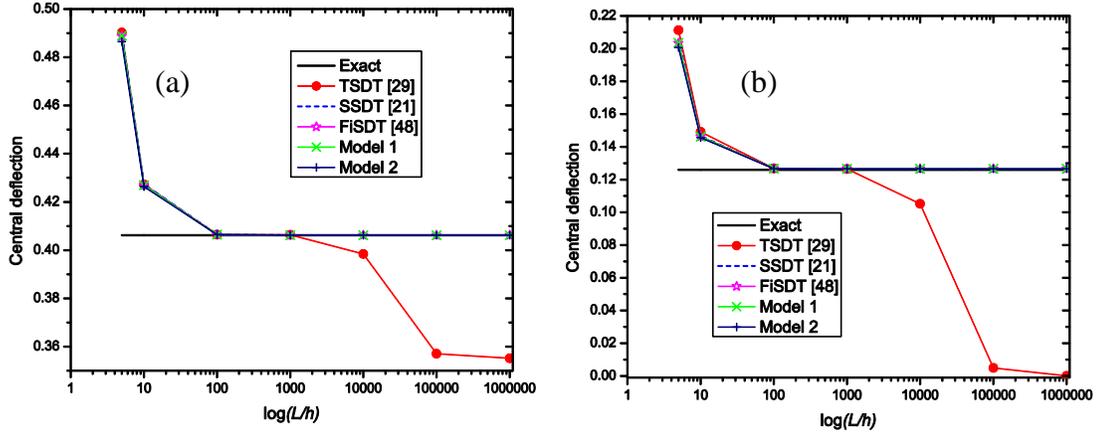

Figure 6. The center deflection via length to thickness ratio under boundary conditions: (a) full support (SSSS), (b) full clamped (CCCC).

*4.2 Static analysis*

In this example, the Al/Al$_2$O$_3$ FGM plate subjected to a sinusoidal pressure defined as $q_0 \sin(\frac{\pi x}{a})\sin(\frac{\pi y}{a})$ is considered. It is noted that the material properties are computed by the rule of mixture based on Eq. (2). Using different functions *f(z)* given in Table 1, the results based on RPT including deflection $\bar{w}$ and axial stress $\bar{\sigma}_x(\frac{h}{3})$ at plate center are summarized in Table 3. The solutions are in good agreement with that of Zenkour's generalized shear deformation theory [52], and those from Carrera et al. [53,58] using Carrera's unified formulation and Neves et al. [54,55] using sinusoidal shear deformation theory (SSDT) and HSDT models. It is concluded that the quasi 3D models accounting for the thickness stretching effect $\varepsilon_z \neq 0$ gain the lower transverse displacement and higher axial stress than the 2D plate models which eliminate the stretching effect. However, the discrepancy between two models reduces as the plate becomes thinner. Figure 7 plots the stress distribution through the thickness of thick plate with *a/h* = 4 and *n* = 1. Using RPT models, the axial stress is plotted in the same path while there is a slight difference observed for shear stress distribution. And, all of them satisfy the traction-free boundary conditions at the plate surfaces.

Table 3: The non-dimensional deflection and, axial stress of SSSS Al/Al$_2$O$_3$ square plate under sinusoidal load.

| *n* | Model | $\varepsilon_z$ | *a/h* = 4 | | 10 | | 100 | |
|---|---|---|---|---|---|---|---|---|
| | | | $\bar{w}$ | $\bar{\sigma}_x(z=\frac{h}{3})$ | $\bar{w}$ | $\bar{\sigma}_x(z=\frac{h}{3})$ | $\bar{w}$ | $\bar{\sigma}_x(z=\frac{h}{3})$ |
| | Ref. [58] | ≠0 | 0.7171 | 0.6221 | 0.5875 | 1.5064 | 0.5625 | 14.969 |
| | CLT | 0 | 0.5623 | 0.806 | 0.5623 | 2.015 | 0.5623 | 20.15 |
| | FSDT | 0 | 0.7291 | 0.806 | 0.5889 | 2.015 | 0.5625 | 20.15 |
| 1 | | GSDT [52] | 0 | - | - | 0.5889 | 1.4894 | - | - |
| | | Ref.[53] | 0 | 0.7289 | 0.7856 | 0.589 | 2.0068 | 0.5625 | 20.149 |
| | | Ref.[53] | ≠0 | 0.7171 | 0.6221 | 0.5875 | 1.5064 | 0.5625 | 14.969 |
| | | SSDT [54] | ≠0 | 0.6997 | 0.5925 | 0.5845 | 1.4945 | 0.5624 | 14.969 |
| | | HSDT [55] | 0 | 0.7308 | 0.5806 | 0.5913 | 1.4874 | 0.5648 | 14.944 |
| | | HSDT [55] | ≠0 | 0.702 | 0.5911 | 0.5868 | 1.4917 | 0.5647 | 14.945 |
| | RPT | TSDT [10] | | 0.7284 | 0.5796 | 0.5889 | 1.4856 | 0.5625 | 14.9255 |
| | | SSDT [21] | | 0.728 | 0.5787 | 0.5889 | 1.4852 | 0.5625 | 14.9255 |
| | | HSDT [20] | | 0.7271 | 0.5779 | 0.5888 | 1.4849 | 0.5625 | 14.9255 |
| | | FiSDT[48] | | 0.725 | 0.5765 | 0.5885 | 1.4844 | 0.5625 | 14.9254 |
| | | Model 1 | | 0.7254 | 0.5779 | 0.5885 | 1.4849 | 0.5625 | 14.9255 |
| | | Model 2 | | 0.7204 | 0.5793 | 0.5878 | 1.4854 | 0.5625 | 14.9255 |
| 4 | | Ref. [58] | ≠0 | 1.1585 | 0.4877 | 0.8821 | 1.1971 | 0.8286 | 11.923 |
| | | CLT | 0 | 0.8281 | 0.642 | 0.8281 | 1.6049 | 0.8281 | 16.049 |
| | | FSDT | 0 | 1.1125 | 0.642 | 0.8736 | 1.6049 | 0.828 | 16.049 |
| | | GSDT [52] | 0 | - | - | 0.8651 | 1.1783 | - | - |
| | | Ref. [53] | 0 | 1.1673 | 0.5986 | 0.8828 | 1.5874 | 0.8286 | 16.047 |
| | | Ref. [53] | ≠0 | 1.1585 | 0.4877 | 0.8821 | 1.1971 | 0.8286 | 11.923 |
| | | SSDT [54] | ≠0 | 1.1178 | 0.4404 | 0.875 | 1.1783 | 0.8286 | 11.932 |
| | | HSDT [55] | 0 | 1.1553 | 0.4338 | 0.877 | 1.1592 | 0.8241 | 11.737 |
| | | HSDT [55] | ≠0 | 1.1108 | 0.433 | 0.87 | 1.1588 | 0.824 | 11.737 |
| | RPT | TSDT [10] | | 1.1599 | 0.4433 | 0.8815 | 1.1753 | 0.8287 | 11.8796 |
| | | SSDT [21] | | 1.1619 | 0.4408 | 0.8819 | 1.1742 | 0.8287 | 11.8796 |
| | | HSDT [20] | | 1.1627 | 0.4385 | 0.8821 | 1.1733 | 0.8287 | 11.8796 |
| | | FiSDT[48] | | 1.1614 | 0.4349 | 0.8819 | 1.1718 | 0.8287 | 11.8792 |
| | | Model 1 | | 1.162 | 0.4371 | 0.882 | 1.1727 | 0.8287 | 11.8793 |
| | | Model 2 | | 1.1562 | 0.4369 | 0.8812 | 1.1726 | 0.8287 | 11.8793 |
| 10 | | Ref. [58] | ≠0 | 1.3745 | 0.3695 | 1.0072 | 0.8965 | 0.9361 | 8.9077 |
| | | CLT | 0 | 0.9354 | 0.4796 | 0.9354 | 1.199 | 0.9354 | 11.99 |
| | | FSDT | 0 | 1.3178 | 0.4796 | 0.9966 | 1.199 | 0.936 | 11.99 |
| | | GSDT [52] | 0 | - | - | 1.0089 | 0.8775 | - | - |
| | | Ref. [53] | 0 | 1.3925 | 0.4345 | 1.009 | 1.1807 | 0.9361 | 11.989 |
| | | Ref. [53] | ≠0 | 1.3745 | 0.1478 | 1.0072 | 0.8965 | 0.9361 | 8.9077 |
| | | SSDT [54] | ≠0 | 1.349 | 0.3227 | 0.875 | 1.1783 | 0.8286 | 11.932 |
| | | HSDT [55] | 0 | 1.376 | 0.3112 | 0.9952 | 0.8468 | 0.9228 | 8.6011 |
| | | HSDT [55] | ≠0 | 1.3334 | 0.3097 | 0.9888 | 0.8462 | 0.9227 | 8.601 |
| | RPT | TSDT [10] | | 1.3908 | 0.3249 | 1.0087 | 0.876 | 0.9362 | 8.8804 |
| | | SSDT [21] | | 1.3917 | 0.3225 | 1.0089 | 0.875 | 0.9362 | 8.8804 |
| | | HSDT[20] | | 1.3906 | 0.3203 | 1.0088 | 0.8741 | 0.9362 | 8.8804 |
| | | FiSDT[48] | | 1.3862 | 0.317 | 1.0083 | 0.8727 | 0.9362 | 8.8801 |
| | | Model 1 | | 1.3871 | 0.3189 | 1.0084 | 0.8735 | 0.9362 | 8.8802 |
| | | Model 2 | | 1.3738 | 0.3183 | 1.0064 | 0.8732 | 0.9362 | 8.8802 |

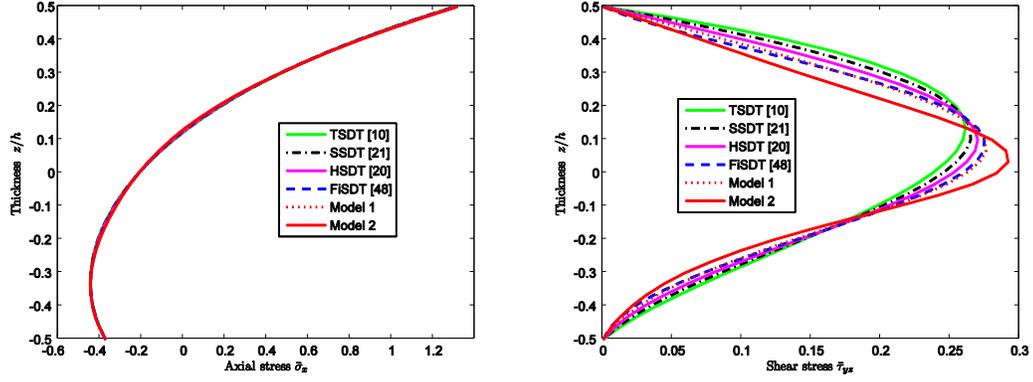

Figure 7. The stresses through thickness of Al/Al2O3 FG plate under sinusoidal load with *a/h*=4, *n*=1, via different refined plate models.

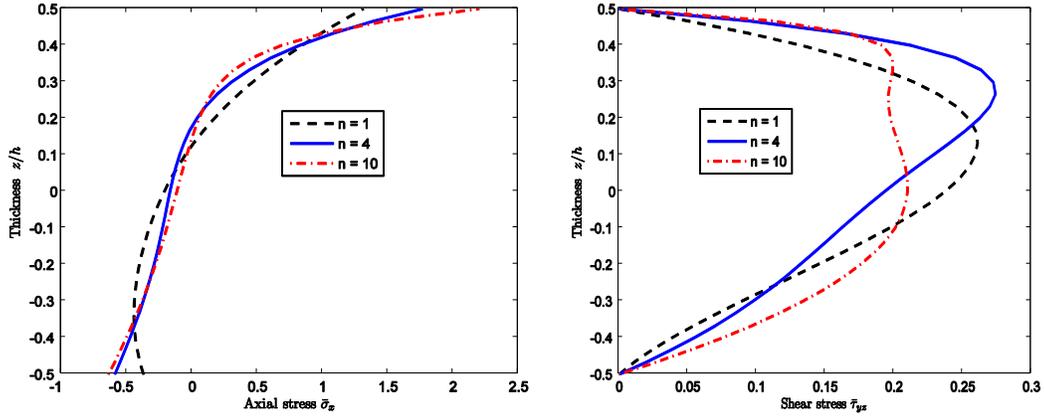

Figure 8. The stresses through thickness of Al/Al2O3 FG plate under sinusoidal load with *a/h*=4, *n*=1 via various power indices *n*.

Figure 8 presents the curved distribution of the axial and shear stresses through the plate thickness according to power index $n = 1, 4, 10$, respectively. It can be concluded that the present model based on NURBS approximation yields very promising results compared to that given in Ref. [55].

Table 4: The non-dimension deflection of Al/ZrO$_2$-1 plate under uniform load with *a/h*=5 via different boundary conditions.

| BC | n | Model | | | | |
|---|---|---|---|---|---|---|
| | | HOSNDPT[59] | TSDT[29] | PRT | | |
| | | | | FiSDT[48] | Model 1 | Model 2 |
| | ceramic | 0.5019 | 0.5088 | 0.5073 | 0.5074 | 0.506 |
| | 0.5 | 0.7543 | 0.7607 | 0.7587 | 0.7588 | 0.7568 |
| | 1 | 0.8708 | 0.8776 | 0.8754 | 0.8756 | 0.8732 |
| SFSF | 2 | 0.9744 | 0.9830 | 0.9808 | 0.981 | 0.9784 |
|      | 4       | -      | 1.0701 | 1.0676 | 1.0679 | 1.0648 |
|      | 8       | -      | 1.1577 | 1.1540 | 1.1544 | 1.1504 |
|      | metal   | 1.4345 | 1.4537 | 1.4490 | 1.4498 | 1.4458 |
|      | ceramic | 0.1671 | 0.1716 | 0.1711 | 0.1711 | 0.1703 |
|      | 0.5     | 0.2505 | 0.2554 | 0.2547 | 0.2548 | 0.2536 |
|      | 1       | 0.2905 | 0.2955 | 0.2947 | 0.2948 | 0.2934 |
| SSSS | 2       | 0.328  | 0.3334 | 0.3327 | 0.3328 | 0.3312 |
|      | 4       | -      | 0.3655 | 0.3647 | 0.3649 | 0.363  |
|      | 8       | -      | 0.3958 | 0.3943 | 0.3945 | 0.3922 |
|      | metal   | 0.4775 | 0.4903 | 0.4887 | 0.4889 | 0.4865 |
|      | ceramic | 0.0731 | 0.0734 | 0.0709 | 0.071  | 0.0701 |
|      | 0.5     | 0.1073 | 0.1077 | 0.1041 | 0.1043 | 0.1029 |
|      | 1       | 0.1253 | 0.1256 | 0.1215 | 0.1217 | 0.1201 |
| CCCC | 2       | 0.1444 | 0.1447 | 0.1401 | 0.1402 | 0.1384 |
|      | 4       | -      | 0.1622 | 0.1566 | 0.1568 | 0.1546 |
|      | 8       | -      | 0.1760 | 0.1694 | 0.1696 | 0.1669 |
|      | metal   | 0.2088 | 0.2098 | 0.2027 | 0.203  | 0.2001 |

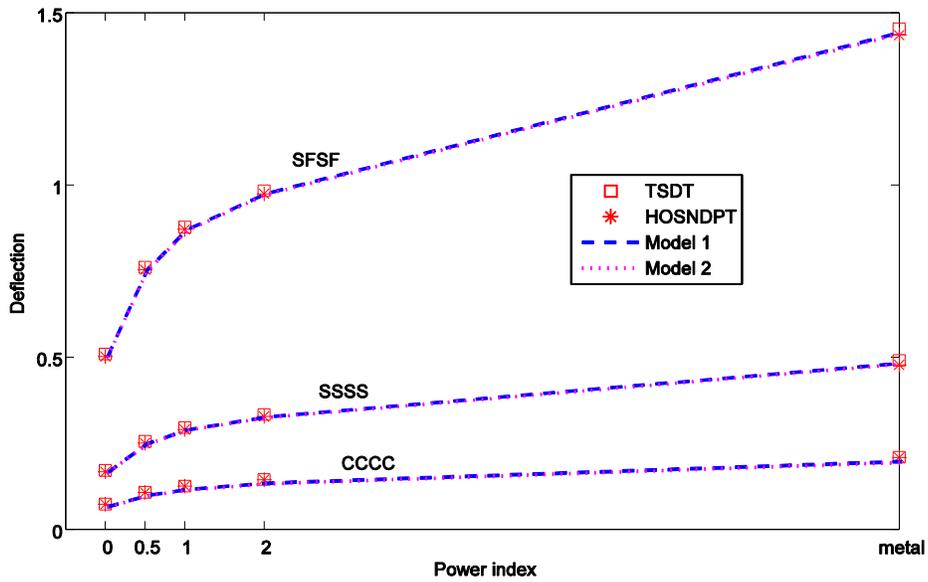

Figure 9. The normalized deflection of Al/ZrO2-1 FGM plate via power indexes and boundary conditions.

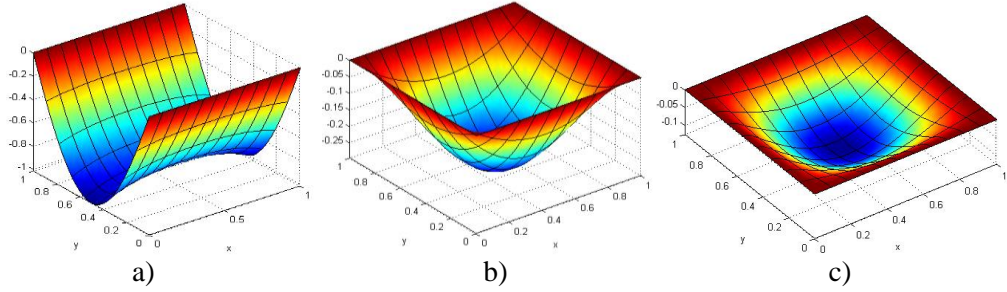

Figure 10. Deflection profile of Al/ZrO$_2$-1 FGM plates: (a) SFSF; (b) SSSS; (c) CCCC.

To end this sub-section, the effective of boundary conditions on the normalized central deflection $\bar{w}_c = 100 w_c E_m h^3 / \{12(1-\nu^2) q_0 a^4\}$ of Al/ZrO$_2$-1 plate is investigated. In this problem the Mori-Tanaka scheme is used for homogenizing Al/ZrO$_2$-1 FGM plate. The present results listed in Table 4 are compared with those of Gilhooley[59] based on HOSNDPT using 18DOFs/node and Tran *et al*. [29] based on TSDT using 5 DOFs/node. Figure 9 shows that, by just using 4 DOFs/node, the present model produces solutions very close to Refs. [29,59] for all boundary conditions. Moreover, when the boundary condition changes from CCCC to SSSS and SFSF, the structural stiffness reduces, the magnitudes of deflection thus increase, respectively. The shapes of transverse displacement according the various boundary conditions are illustrated in Figure 10.

*4.3 Free vibration analysis*

Let us consider a simply supported Al/ZrO$_2$-1 plate which is homogenized by the Mori-Tanaka scheme. In Table 5, we also provide the results based on the RPT with various functions $f(z)$ in [10, 48] and two current models 1 & 2 in comparison with the exact solution [50], that of HOSNDPT [47] and quasi 3D-solution using SSDT and HSDT [54, 55]. The excellent correlation between these models is again achieved for all values of exponent *n*. In addition, it is revealed that the proposed model 2 gives the best natural frequency with the least error compared to exact result [50] by Vel and Batra. The first ten natural frequencies of the thick and moderate plate with *a/h* = 5, 10, 20 are listed in Table 6. The computed values agree well with the literature [47] for various *a/h* ratios and mode number. Corresponding to *n* = 1, the first six mode shapes are plotted in Figure 11.

Table 5: The natural frequency $\bar{\omega} = \omega h \sqrt{\rho_m / E_m}$ of SSSS Al/ZrO$_2$-1 plate with *a/h*=5.

| Model | $\varepsilon_z$ | \multicolumn{7}{c}{n} |
|---|---|---|---|---|---|---|---|---|
| | | 0 | 0.5 | 1 | 2 | 3 | 5 | 10 |
| Exact [50] | - | - | - | 0.2192 | 0.2197 | 0.2211 | 0.2225 | - |
| HOSNDPT [47] | - | - | - | 0.2152 (-1.82)[*] | 0.2153 (-2.00) | 0.2172 (-1.76) | 0.2194 (-1.39) | - |

|  |  |  |  |  | 0.2184 | 0.2189 | 0.2202 | 0.2215 |  |
|---|---|---|---|---|---|---|---|---|---|
| SSDT [54] |  | 0 | - | - | (-0.36) | (-0.36) | (-0.41) | (-0.45) | - |
| SSDT [54] |  | ≠0 | - | - | 0.2193 | 0.2198 | 0.2212 | 0.2225 | - |
|  |  |  |  |  | (0.05) | (0.05) | (0.05) | (0.00) |  |
| HSDT [55] |  | 0 | 0.2459 | 0.2219 | 0.2184 | 0.2191 | 0.2206 | 0.222 | 0.2219 |
|  |  |  |  |  | (-0.36) | (-0.27) | (-0.23) | (-0.22) |  |
| HSDT [55] |  | ≠0 | 0.2469 | 0.2228 | 0.2193 | 0.22 | 0.2215 | 0.223 | 0.2229 |
|  |  |  |  |  | (0.05) | (0.14) | 90.18) | (0.22) |  |
| RPT | TSDT[10] |  | 0.2459 | 0.2221 | 0.2184 | 0.2189 | 0.2203 | 0.2216 | 0.2211 |
|  |  |  |  |  | (-0.36) | (-0.36) | (-0.36) | (-0.40) |  |
|  | FiSDT[48] |  | 0.2462 | 0.2224 | 0.2187 | 0.2191 | 0.2205 | 0.2218 | 0.2215 |
|  |  |  |  |  | (-0.23) | (-0.27) | (-0.27) | (-0.31) |  |
|  | Model 1 |  | 0.2462 | 0.2224 | 0.2186 | 0.2191 | 0.2205 | 0.2218 | 0.2215 |
|  |  |  |  |  | (-0.27) | (-0.27) | (-0.27) | (-0.31) |  |
|  | Model 2 |  | 0.2468 | 0.2229 | 0.2192 | 0.2196 | 0.221 | 0.2224 | 0.2222 |
|  |  |  |  |  | (0.00) | (-0.05) | (-0.05) | (-0.04) |  |

(*) The error in parentheses

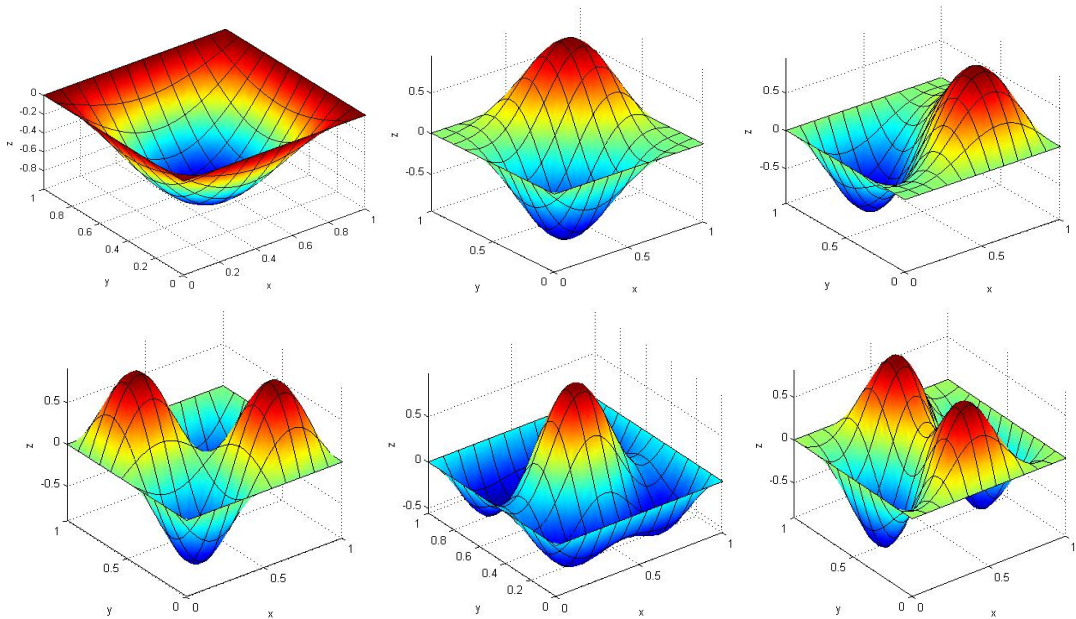

Figure 11 The first six mode shapes of Al/ZrO$_2$-1 with $n = 1$, $a/h=5$.

Table 6: The natural frequency $\bar{\omega}$ of SSSS Al/ZrO$_2$-1 plate with various ratios $a/h$.

| $a/h$ | Model |  | Modes | | | | | | | |
|---|---|---|---|---|---|---|---|---|---|---|
|  |  | 1 | (2,3) | 4 | 5 | 6 | 7 | 8 | 9 | 10 |
| 5 | Exact [50] | 0.2192 | - | - | - | - | - | - | - | - |
|  | HOSNDPT [47] | 0.2152 | 0.4114 | 0.4761 | 0.4761 | 0.582 | 0.6914 | 0.8192 | 0.8217 | 0.8242 |
|  | FiSDT[48] | 0.2187 | 0.4116 | 0.4806 | 0.4806 | 0.5821 | 0.6976 | 0.8233 | 0.8233 | 0.8263 |
| RPT | Model 1 | 0.2186 | 0.4116 | 0.4804 | 0.4804 | 0.5821 | 0.6972 | 0.8233 | 0.8233 | 0.8257 |
|    | Model 2     | 0.2192 | 0.4116 | 0.4827 | 0.4827 | 0.5821 | 0.7018 | 0.8233 | 0.8233 | 0.832  |
| 10 | Exact [50]  | 0.0596 | -      | -      | -      | -      | -      | -      | -      | -      |
|    | HOSNDPT [47]| 0.0584 | 0.141  | 0.2058 | 0.2058 | 0.2164 | 0.2646 | 0.2677 | 0.2913 | 0.3264 |
|    | SSDT [54]   | 0.0596 | 0.1426 | 0.2058 | 0.2058 | 0.2193 | 0.2676 | 0.2676 | 0.291  | 0.3363 |
|    | HSDT [55]   | 0.0596 | 0.1426 | 0.2059 | 0.2059 | 0.2193 | 0.2676 | 0.2676 | 0.2912 | 0.3364 |
|    | FiSDT [48]  | 0.0595 | 0.1423 | 0.2058 | 0.2058 | 0.2187 | 0.2668 | 0.2668 | 0.2911 | 0.3351 |
| RPT| Model 1     | 0.0595 | 0.1423 | 0.2058 | 0.2058 | 0.2187 | 0.2667 | 0.2667 | 0.2911 | 0.335  |
|    | Model 2     | 0.0596 | 0.1425 | 0.2058 | 0.2058 | 0.2192 | 0.2675 | 0.2675 | 0.2911 | 0.3362 |
| 20 | Exact [28]  | 0.0153 | -      | -      | -      | -      | -      | -      | -      | -      |
|    | HOSNDPT [47]| 0.0149 | 0.0377 | 0.0593 | 0.0747 | 0.0747 | 0.0769 | 0.0912 | 0.0913 | 0.1029 |
|    | SSDT [54]   | 0.0153 | 0.0377 | 0.0596 | 0.0739 | 0.0739 | 0.095  | 0.095  | 0.1029 | 0.1029 |
|    | HSDT [55]   | 0.0153 | 0.0377 | 0.0596 | 0.0739 | 0.0739 | 0.095  | 0.095  | 0.103  | 0.103  |
|    | FiSDT [48]  | 0.0153 | 0.0377 | 0.0595 | 0.0739 | 0.0739 | 0.0949 | 0.0949 | 0.1029 | 0.1029 |
| RPT| Model 1     | 0.0153 | 0.0377 | 0.0595 | 0.0739 | 0.0739 | 0.0949 | 0.0949 | 0.1029 | 0.1029 |
|    | Model 2     | 0.0153 | 0.0377 | 0.0596 | 0.0739 | 0.0739 | 0.095  | 0.095  | 0.1029 | 0.1029 |

*4.4 Buckling analysis*

In this section, a clamped circular plate of radius $R$ and thickness $h$ is subjected to a uniform radial pressure $p_0$ is meshed into 11x11 cubic elements as shown in Figure 12. The plate made from Al/ZrO$_2$-2 is homogenized following the rule of mixture, in which the effective Young's modulus $E_e$ and Poison's ratio $\nu_e$ are calculated according to Eq. (35).

$$P_e = P_c V_c(z) + P_m V_m(z) \text{ where } V_m(z) = \left(\frac{1}{2} - \frac{z}{h}\right)^n, \quad V_c = 1 - V_m \tag{35}$$

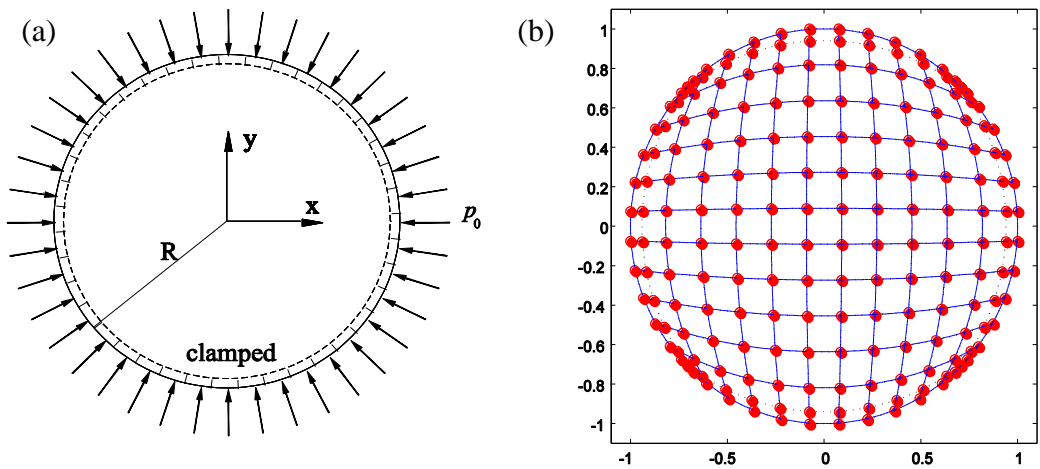

Figure 12. The circular plate: (a) geometry and (b) meshing of 11x11 cubic elements.

The variation of the critical load parameter $\bar{p}_{cr}=p_{cr}R^2/D_m$ versus the exponential value and the thickness to radius ratio is illustrated in Figure 13 and Figure 14 respectively. It is observed that the critical buckling load increases corresponding to the decrease in the *h/R* ratios and the increase in the value of *n*. When *n>* 10, it is slightly independent on the increment of the exponent *n*. Table 7 shows the comparison between present results with analytical solution published in the literatures [56,57] and numerical solution based on IGA [29] using TSDT. The RPT gives slightly higher values than TSDT [57,29] and unconstrained third-order shear deformation plate theory (UTSDT) [56] with the difference of 0.5% and 1.8%. However, the discrepancy between them reduces as plate becomes thinner. Figure 15 presents the first four buckling mode shapes of the clamped circular Al/ZrO$_2$-2 plate with *h/R* = 0.1 and *n* = 2.

Table 7: Comparison of the buckling load parameter of clamped thick circular Al/ZrO$_2$-2 plate.

| Power index | Method | | *h/R* | | | |
|---|---|---|---|---|---|---|
| | | | 0.1 | 0.2 | 0.25 | 0.3 |
| 0 | FST [56] | | 14.089 | 12.571 | 11.631 | 10.657 |
| | TST [57] | | 14.089 | 12.574 | 11.638 | 10.67 |
| | UTSDT[56] | | 14.089 | 12.575 | 11.639 | 10.67 |
| | TSDT[29] | | 14.1089 | 12.5914 | 11.6540 | 10.6842 |
| | RPT | FiSDT[48] | 14.1873 | 12.6787 | 11.7466 | 10.7822 |
| | | Model 1 | 14.1859 | 12.6743 | 11.7405 | 10.7745 |
| | | Model 2 | 14.2023 | 12.7281 | 11.8143 | 10.8666 |
| 0.5 | FST[56] | | 19.423 | 17.34 | 16.048 | 14.711 |
| | TST[57] | | 19.411 | 17.311 | 16.013 | 14.672 |
| | UTSDT[56] | | 19.413 | 17.31 | 16.012 | 14.672 |
| | TSDT[29] | | 19.4391 | 17.3327 | 16.0334 | 14.6910 |
| | | FiSDT[48] | 19.5458 | 17.4504 | 16.1579 | 14.8227 |
| | RPT | Model 1 | 19.5439 | 17.4441 | 16.1492 | 14.8118 |
| | | Model 2 | 19.5663 | 17.518 | 16.2506 | 14.9381 |
| 2 | FST[56] | | 23.057 | 20.742 | 19.29 | 17.77 |
| | TST[57] | | 23.074 | 20.803 | 19.377 | 17.882 |
| | UTSDT[56] | | 23.075 | 20.805 | 19.378 | 17.881 |
| | IGA-TSDT[29] | | 23.1062 | 20.8319 | 19.4033 | 17.9060 |
| | | FiSDT[48] | 23.2361 | 20.9794 | 19.5612 | 18.0745 |
| | RPT | Model 1 | 23.2342 | 20.9728 | 19.552 | 18.0628 |
| | | Model 2 | 23.2592 | 21.0569 | 19.6687 | 18.2099 |
| 5 | FST[56] | | 25.411 | 22.876 | 21.282 | 19.611 |
| | TST[57] | | 25.439 | 22.971 | 21.414 | 19.78 |
| | UTSDT[56] | | 25.442 | 22.969 | 21.412 | 19.778 |
| | IGA-TSDT[29] | | 25.4743 | 22.99918 | 21.4407 | 19.80426 |
| | | FiSDT[48] | 25.6172 | 23.1598 | 21.6118 | 19.9861 |
| | RPT | Model 1 | 25.6152 | 23.1529 | 21.6022 | 19.9738 |
| | | Model 2 | 25.6418 | 23.2426 | 21.7268 | 20.1313 |

| | | | | | |
|---|---|---|---|---|---|
| | FST[56] | 27.111 | 24.353 | 22.627 | 20.823 |
| | TST[57] | 27.133 | 24.423 | 22.725 | 20.948 |
| | UTSDT[56] | 27.131 | 24.422 | 22.725 | 20.949 |
| 10 | TSDT[29] | 27.1684 | 24.4542 | 22.7536 | 20.9750 |
| | FiSDT[48] | 27.3176 | 24.6148 | 22.9216 | 21.1509 |
| PRT | Model 1 | 27.3155 | 24.6077 | 22.9117 | 21.1383 |
| | Model 2 | 27.3429 | 24.6994 | 23.0389 | 21.2986 |

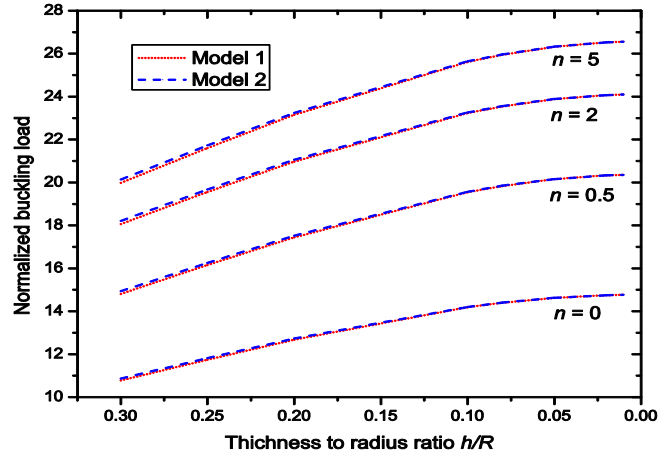

Figure 13. The normalized buckling load $\bar{p}_{cr}$ of the Al/ZrO$_2$-2 plate via radius to thickness ratio $R/h$.

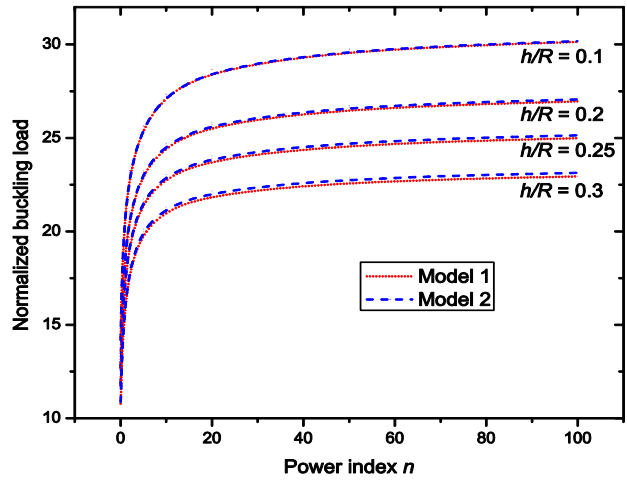

Figure 14. The normalized buckling load $\bar{p}_{cr}$ of the Al/ZrO$_2$-2 plate via power index $n$.

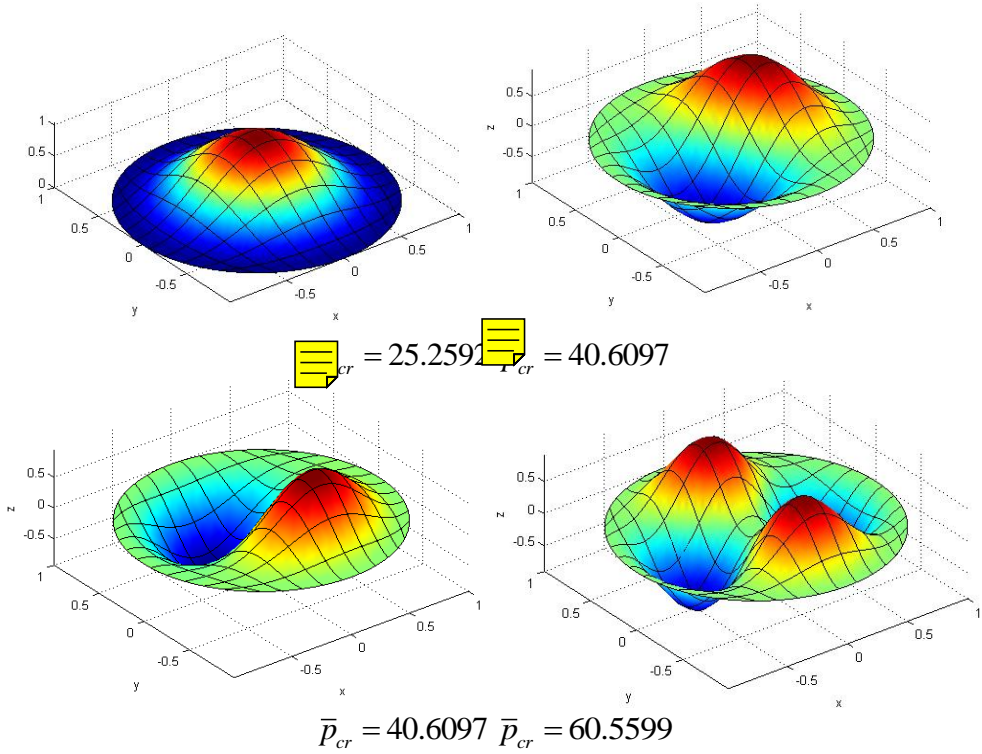

$$\bar{p}_{cr} = 25.2592 \quad \bar{p}_{cr} = 40.6097$$

$$\bar{p}_{cr} = 40.6097 \quad \bar{p}_{cr} = 60.5599$$

Figure 15. The first four buckling modes of Al/ZrO$_2$-2 plate with $h/R = 0.1$ and $n = 2$.

## 5. Conclusions

In this paper, we have developed the inverse tangent transverse shear deformation model together with the isogeometric finite element formulation for static, free vibration and buckling analyses of FGM plates. The method fulfills the $C^1$ – requirement of RPT model and the approximate displacements have four DOFs per each control point which results in less computational cost compared to other five DOFs methods. In the present model, the shear strains/stresses are obtained independently on the bending component. As a result, RPT has been strongly similar to CPT and it naturally overcomes the shear locking phenomenon. The present results are compared with analytical solutions and those using HSDT or quasi-3D models and they demonstrated excellent agreement in static, free vibration and buckling problems. The present models ensured the non-linear distribution of the shear stress/strain through the plate thickness without using any shear correction factors, and yielded the traction-free boundary conditions at plate surfaces. Finally, it can be concluded that the proposed model 2 produces the high accuracy for analysis of FGM plates.